\documentclass[11pt]{article}

\usepackage{latexsym,amsfonts,amssymb,amsmath,amsthm}
\usepackage{graphicx}
\usepackage{color}
\usepackage{mathrsfs}
\usepackage{amsmath, amsthm, amsfonts, amssymb, cite}
\usepackage{wrapfig}
\usepackage{stackrel}
\usepackage{color}
\usepackage{float}
\usepackage{enumitem}
\usepackage{mathtools}
\usepackage{multicol}
\usepackage[normalem]{ulem}
\usepackage{amsmath,amsthm,verbatim,amssymb,amsfonts,amscd,graphicx}
\usepackage{graphics}
\usepackage{relsize}
\usepackage{esvect}
\usepackage{mathtools}
\usepackage{physics}
\usepackage{multirow}
\usepackage{hyperref}
\usepackage{xcolor}

\theoremstyle{plain}
\newtheorem{theorem}{Theorem}
\newtheorem*{theorem*}{Theorem}
\newtheorem{corollary}{Corollary}
\newtheorem{lemma}{Lemma}

\theoremstyle{definition}

\newtheorem{remark}{Remark}

\begin{document}

\title{On the regularized $L^4$-norm for Eisenstein series in the level aspect, Part II}
\author{Jiakun Pan}
\maketitle

\section{Introduction}
This is a sequel paper of \href{https://arxiv.org/abs/2003.10995}{\underline{\textit{On the regularized $L^4$-norm for Eisenstein series in the level aspect}}} by the same author. The two papers will be combined for publication.
\subsection{Background}
The Gaussian Moments Conjecture is a number theoretical manifestation of the Random Wave Conjecture on the randomness of automorphic forms in the spectral aspect. For Eisenstein series, we can either consider the regularized moments of them, or the moments of its truncation. On Eisenstein series of fixed levels, we have the results on the fourth moments by Spinu\cite{Sp}, Humphries\cite{H}, Djankovi\'c and Khan\cite{DK1,DK2}.
In the previous paper, we let the level to increase and reduced the fourth moment of newform Eisenstein series attached to a cusp with
\begin{align*}
    \langle |E_{\mathfrak{a}}(\cdot, \frac{1}{2}+iT, \chi)|^2, |E_{\mathfrak{a}}(\cdot, \frac{1}{2}+iT, \chi)|^2 \rangle_{\mathrm{reg}} = I_1 + I_2,
\end{align*}
where
\begin{align*}
    I_{1}= \nu^{-1}(N) \sum_{u}  \frac{\Lambda^2(\frac{1}{2},u) |\Lambda(\frac{1}{2}+2iT, u\otimes \psi)|^2}{|\Lambda(1+2iT,\psi)|^2} + \mathrm{continuous} \medspace\medspace  \mathrm{spectrum},
\end{align*}
for some primitive $\psi$ mod $N$ delicately decided by $\chi$ and $\mathfrak{a}$, and
\begin{align}\label{part1}
    I_2 = \nu^{-1}(N) \Big( \frac{24}{\pi}\log^2 N + O(\frac{L''}{L}(1+2iT,\psi)) + O(\log N \log\log N \frac{L'}{L}(1+2iT, \psi)) \Big).
\end{align}
Here we multiply $I_1$ and $I_2$ by $\nu(N):=[SL_2(\mathbb{Z}):\Gamma_0(N)]$ for clarity, and the sum over $u$ traverses an orthonormal basis of the discrete spectrum spanned by Maa{\ss} forms.

In this article, we present a more hypothetical approach to see what the main terms shall be. Before we state the results, we first describe the conditionality of our research.

\subsection{Assumptions}
For newform Eisenstein series, the period integrals ofteh involve logarithmic derivatives of Dirichlet L-functions. When studying the Quantum Unique Ergodicity (\textbf{QUE}) problem for newform Eisenstein series, the author and Young \cite[Thm 1.3]{PY} found that
\begin{align*}
    \frac{\langle |E|^2, \phi \rangle}{\langle 1, \phi \rangle} = \frac{6}{\pi}\log N + \frac{12}{\pi}  \mathrm{Re}\frac{L'}{L}(1+2iT,\psi) + O(N^{-\frac{1}{8}+\varepsilon}),
\end{align*}
where $E=E_{\mathfrak{a}}(z,\tfrac{1}{2}+iT,\chi)$ for primitive $\chi$ mod $N$, $\psi=\psi(\chi,\mathfrak{a})$ is also primitive mod $N$, and $\phi$ is an $SL_2(\mathbb{Z})$-automorphic test function with nice analytic properties. An unconditional $o(\log N)$ bound being unavailable for the moment, the Generalized Riemann Hypothesis (\textbf{GRH}) implies the following inequalities, which is also useful in this paper:
\begin{align}\label{A1}
    \begin{split}
        \frac{L'}{L}(1+2iT,\chi) &\ll \log\log N; \\
        \frac{L''}{L}(1+2iT,\chi) &\ll (\log\log N)^2.
    \end{split}
\end{align}
GRH also implies the Linderl\"of Hypothesis, which helps to bound some L-functions.

In addition, our main theorem heavily relies on a recipe of conjecturing the moments of L-functions, due to Conrey, Farmer, Keating, Rubinstein, and Snaith \cite{CFKRS}, see Section \ref{recipe} for a sketch. 

\subsection{Statement of main result}
Now we come back to the spectral sum $I_1$ obtained from the regularized fourth moment of weight $0$ newform Eisenstein series $E=E_{\mathfrak{a}}(z,\tfrac{1}{2}+iT,\chi)$, for fixed $T\in \mathbb{R}$.  Recall that $\chi$ is even, primitive mod $N$, and $\mathfrak{a}$ is singular for $\chi$.
\begin{theorem}\label{main}
Assume GRH and the recipe in \cite{CFKRS}. For all prime $N$ we have
\begin{align*}
    I_1 \sim \frac{24}{\pi} \frac{\log^2 N}{\nu(N)} \Big(1 + \delta_{\chi^2=\chi_{0,N}} \delta_{T=0} \Big).
\end{align*}
\end{theorem}
\noindent Combining with our unconditional estimation for $I_2$, we have the following corollary.
\begin{remark}\label{composite}
Indeed we can relax $N$ to all positive integers with
\begin{align*}
    \sigma_{-1}(N) = 1+O(N^{-\delta})
\end{align*}
for any fixed $\delta>0$. See Remark \ref{relaxation} for an explanation.
\end{remark}

\begin{remark}\label{varyTQ}
Again, we have assumed $T$ to be fixed here. Since each Eisenstein series $E$ is continuous in $T$, it makes us curious about what happens when $T\rightarrow 0$ as $N$ grows. For this sake we allow $T=T(N)$ to approach zero, and conclude that $T\asymp \log^{-1} N$ makes the threshold. In other words, depending on whether $T$ shrinks faster or slower than $\log^{-1} N$, $I_1=I_1(T)$ has different asymptotics.
\end{remark}

\begin{corollary}\label{rvc}
With the same assumptions as in Theorem \ref{main}, we have
\begin{align}\label{2cases}
    \begin{split}
        \frac{\langle |E|^2, |E|^2 \rangle_{\mathrm{reg}}}{(\mathrm{Vol}(\Gamma_0(N)\backslash\mathbb{H}))^{-1} (\sqrt{2\log N})^4} \sim     2 \cdot \begin{cases} 
 3 & \text{when } \chi \text{ is quadratic and } T=0; \\ 2 & \text{otherwise}.
    \end{cases}
    \end{split}
\end{align}
\end{corollary}

The left hand side can be regarded as the fourth moment of $E$ under proper rescaling; the two cases on the right hand side correspond to two models of the Gaussian Moments Conjecture.
When $T=0$ and $\chi$ is quadratic, there exists a complex scalar $\epsilon$, such that $\epsilon E$ is real-valued, which is similar to classical Eisenstein series \cite{DK1,DK2} and dihedral Maa{\ss} forms \cite{HK} in the $t$-aspect. So, we expect their moments to behave like a real random wave in the $N$-aspect. Thus, we have a good reason to expect some similarity with the real Gaussian distribution, whose fourth moment is $3$.

In all other cases, the Hecke sequence of $E$ is not contained in any straight line in the complex plane, and as $N$ grows, the limiting behavior resembles the complex Gaussian distribution with fourth moment $2$, as Blomer, Khan and Young showed for holomorphic forms of large weight in \cite{BKY}.

\subsection{The fourth moment of truncated newform Eisenstein series}
Djankovi\'c and Khan \cite{DK1} did a consistency check on their conjecture (later proved in \cite{DK2}) on the regularized fourth moment for classical Eisenstein series, with the Random Wave Conjecture. Following their methods, we obtain an essentially different result in the level aspect. Recall we write $E^Y$ for $E$ truncated in its $Y$-cuspidal zones.

\begin{remark}
Here, we need an analogue of Spinu's \cite{Sp} work\footnote{Up to a couple of errors that are fixed by Humphries \cite{H}} on the optimal upper bound for the fourth moment of truncated classical Eisenstein series, in the level aspect:
\begin{align}\label{A3}
    \langle |E^Y|^2, |E^Y|^2 \rangle \ll_Y \frac{\log^2 N}{\nu(N)}.
\end{align}
Here $E^Y$ stands for the truncated Eisenstein series for $Y>1$; in the spectral parameter aspect, it has been proven by Humphries that the fourth moment of the truncated classical Eisenstein series is $\Omega(\log^2 T)$; GMC says the correct main term should be $\tfrac{36}{\pi}\log^2 T$.

\end{remark}

\begin{theorem}\label{0diff}
Fix $T=0$ and let $\chi$ be quadratic mod $N$.\footnote{The settings on $T$ and $\chi$ is for brevity of discussion; as we point out in Remark \ref{explanation}, things also hold, and are much simpler in other cases.} With GRH and inequality \eqref{A3}, we have
\begin{align*}
    \langle |E|^2, |E|^2 \rangle_{\mathrm{reg}} - \langle |E^Y|^2, |E^Y|^2 \rangle = o_Y(\frac{\log^2 N}{\nu(N)}).
\end{align*}
\end{theorem}

\begin{remark}
A level aspect variant for the Gaussian Moments Conjecture has not been formulated yet, so a first guess would be to simply replace every $t$ with $N$ and renormalize with $\nu^{-1}(N)$ in the results for classical Eisenstein series. However, as Theorem \ref{0diff} suggests, if we expect
\begin{align*}
    \langle |E^Y|^2, |E^Y|^2 \rangle \sim  C_4 \frac{\log^2 N}{\nu(N)},
\end{align*}
for some constant $C_4$ (note this implies \eqref{A3}), then Corollary \ref{rvc} says for $T=0$ and $\chi^2=\chi_{0,N}$,
\begin{align*}
    C_4 = \frac{72}{\pi},
\end{align*}
instead of the na\"ively expected $\tfrac{36}{\pi}$. Why there is an extra factor of $2$ remains an open question.
\end{remark}

\subsection{Acknowledgement}
The author sincerely thanks Max Planck Institute for Mathematics for hosting him to do the research for this paper, as well as his mentor Valentin Blomer for his insightful guidance. He is grateful to Peter Humphries for helpful discussions in the early stage. Finally, he is indebted to his Ph.D. advisor Matthew P. Young for multiple constructive suggestions.

\section{Prerequisite}
\subsection{Automorphic L-functions}
Let $u=u(z)$ be an $L^2$-normalized Maa{\ss} form of trivial nebentypus, spectral parameter $t=t(u)$, and level $M=M(u)$, where $M\mid N$. 
For each $u$, we have the following formal approximate functional equations
    \begin{align}\label{FE1}
    L(s,u) = \sum_{n\geq 1} \frac{\lambda_u(n)}{n^s} g_1(\frac{n}{\sqrt{M}})+ \gamma(u,s) \sum_{n\geq 1} \frac{\lambda_u(n)}{n^{1-s}} g_1(\frac{n}{\sqrt{M}}),
\end{align}
where $\lambda_u(n)$ is the $n$-th Hecke eigenvalue, $g_1$ is some smooth and compactly supported weight function, and $\rho_u \sqrt{y}\lambda_u(n) K_{it}(2\pi |n| y)$ is the $n$-th Fourier coefficient of $u$. Similarly,
\begin{align}\label{FE2}
    L(s,u \otimes \chi) = \sum_{n\geq 1} \frac{\lambda_u(n) \chi(n)}{n^s} g_2(\frac{n}{N}) + \gamma(s, u\otimes \chi) \sum_{n\geq 1} \frac{\lambda_u(n) \overline{\chi}(n)}{n^{1-s}} g_2(\frac{n}{N}),
\end{align}
for some weight function $g_2$.
The $\gamma$-functions are defined as below (see \cite[Sec. 2.1-2.2]{BFKMM}):
\begin{align*}
    \gamma(s,u) &=\lambda_u(-1) f(s,u)\\
    \gamma(s, u\otimes \chi) &= \lambda_u(-1) (\frac{\tau(\chi)}{\sqrt{N}})^2 f(s,u\otimes \chi),
\end{align*}
with (here $Q(u)=M$ and $Q(u\otimes \chi)=N^2$)
\begin{align*}
    f(s,g) &= Q(g)^{\frac{1}{2}-s} \pi^{2s-1} \frac{\Gamma(\frac{1-s+it}{2}) \Gamma(\frac{1-s-it}{2})}{\Gamma(\frac{s+it}{2}) \Gamma(\frac{s-it}{2})}.
\end{align*}
\subsection{Eisenstein series}
The Fourier expression for newform Eisenstein series can be written our explicitly, see \cite{Young} for more details. When $E=E_{\mathfrak{a}}(z,s,\chi)$ with primitive $\chi$ and cusp $\mathfrak{a} \sim_{\Gamma_0(N)} \tfrac{1}{f}$ singular for $\chi$ for some $f\mid N$,\footnote{Here $f$ is uniquely determined by $\mathfrak{a}$, see \cite[Sec 2]{PY} for more details.} with $(f,\tfrac{N}{f})=1$. Then we can uniquely decompose $\chi=\chi_1\overline{\chi_2}$ with $\chi_1, \chi_2$ primitive mod $\tfrac{N}{f}, f$ respectively. Writing $\psi=\chi_1\chi_2$, we have
\begin{align}\label{theta}
    E =\frac{\tau(\chi_2)}{q_2^{-s}} \Lambda^{-1}(2s,\psi) E^*_{\chi_1,\chi_2}(z,s).
\end{align}
The completed Eisenstein series $E^*_{\chi_1,\chi_2}(z,s)$ has the following Fourier expansion (see \cite[Prop 4.1]{Young})
\begin{align}\label{Fourier}
E^* _{\chi_1, \chi_2}(z,s) = e_{\chi_1, \chi_2}^*(y,s) + 2\sqrt{y} \sum_{n \neq 0} \lambda_{\chi_1, \chi_2}(n,s) e(nx) K_{s-\frac{1}{2}}(2\pi |n| y),
\end{align}
where the constant term is
\begin{align*}
e^*_{\chi_1, \chi_2}(y,s)= \delta_{q_1 =1} \theta_{1,\chi_2}(s) (q_2y)^s + \delta_{q_2 =1} \theta_{1,\overline{\chi_1}}(1-s) (q_1y)^{1-s},
\end{align*}
$\lambda_{\chi_1, \chi_2}(n,s) = \chi_2(\frac{n}{|n|}) \sum_{ab=|n|}\chi_1(a)\overline{\chi_2}(b) (\frac{b}{a})^{s-\frac{1}{2}}$, $\tau(\chi)$ is the Gauss sum of $\chi$, and $K_{\alpha}$ is the $K$-Bessel function of order $\alpha \in \mathbb{C}$.
Following this, for any $u$ we have (\cite[Lemma 6.1]{PY})
\begin{align}\label{expression}
    |\langle |E|^2, u \rangle|^2 = |\rho_u|^2 \frac{1+\lambda_u(-1)}{8 N} \frac{|\Gamma(\frac{\frac{1}{2} + it}{2})|^4 \prod_{\pm} |\Gamma(\frac{\frac{1}{2}+2iT \pm it}{2})|^2}{|\Gamma(\frac{1}{2}+iT)|^4} \frac{L^2(\frac{1}{2},u) |L(\frac{1}{2}+2iT, u\otimes \psi)|^2}{|L(1+2iT, \psi)|^4}.
\end{align}
A useful fact from similar calculation is for any $d\mid \tfrac{N}{M}$ we have
\begin{align}\label{oldclass}
    \langle |E|^2, u|_d \rangle = \langle |E|^2, u \rangle.
\end{align}
\subsection{An orthonormal basis}\label{orthobase}
Write $\mathcal{B}_{it}(M)=\mathcal{B}^{(N)}_{it}(M)$ for a collection of Maa{\ss}  newform of level $M$, spectral parameter $it$, and $L^2$-norm $1$ on $\Gamma_0(N)\backslash\mathbb{H}$. Also write $\mathcal{B}(M)=\sqcup_t \mathcal{B}_{it}(M)$. An orthonormal basis of Maa{\ss} forms of level $N$ can be chosen as below:
\begin{align}\label{BM}
    \mathcal{O}(N) := \Big\{ u^{\scriptscriptstyle{<\ell>}}(z)=\sum_{d\mid \ell}\xi_{\ell}(d)u|_d \quad \Big| \quad u \in \mathcal{B}_{it}(M), \ell \mid L, ML=N, t\in \mathrm{Spec}_{\Gamma_0(N)}(\Delta) \Big\},
\end{align}
with $\xi_{\ell}(d)$ defined in Lemma 2.1 of \cite{BM}. In addition we have 
\begin{align}\label{BMsharp}
    \sum_{d\mid \ell}|\xi_{\ell}(d)|=O(\ell^{\varepsilon} ).
\end{align}

\begin{remark}
As Young verified later, the choice of coefficients also makes an orthonormal basis for Eisenstein series, in the sense of the formal inner product $\langle E_{\mathfrak{a}}, E_{\mathfrak{a}}\rangle := \delta_{\mathfrak{a}=\mathfrak{b}} 4\pi$. Henceforth we can let $\mathcal{O}(N)$ contain Eisenstein series and omit specific discussions on them. 
\end{remark}
\subsection{The recipe of conjecturing moments of L-functions}\label{recipe}
There is a recipe to make conjectures on the moments of L-functions, due to B. Conrey, D. Farmer, J. Keating, M. Rubinstein and N. Snaith \cite{CFKRS}. We call it ``the recipe'' for short throughout this article.
To summarize, we can approximate the average of L-functions by
\begin{align*}
    \sum_u \Big|\sum_{n\geq 1} \frac{a_u(n)}{n^s} + \gamma_u \sum_{n\geq 1} \frac{b_u(n)}{n^{1-s}}\Big|^{2k},
\end{align*}
and after we expand the $2k$-th power as
\begin{align*}
     \sum_u \underbrace{\text{Product of Root Numbers}}_{A_u} \cdot \underbrace{\text{Weighted Multi-Dirichlet Series}}_{B_u},
\end{align*}
we can instead estimate
\begin{align*}
    \sum_u (\text{expected mean of } A_u \text{ over the family of } u) \cdot B_u.
\end{align*}
For each type (unitary, orthogonal, symplectic) of L-functions, there are cancellations underneath such that no big loss is expected in the transformation. Needless to say, a rigorous proof that justifies the two expressions to have the same size appears unavailable for the moment.
\begin{remark}
We introduce the recipe for moments of L-functions for convenience, but it is also applicable to hybrid moments of $L$-functions, which is the case for our Theorem \ref{main} with \eqref{FE1} and \eqref{FE2}.
\end{remark}

\subsection{The de Branges-Wilson Beta Integral and the Kuznetsov Trace Formula}
There is an integral formula by Louis de Branges and James A. Wilson, for which we refer the readers to lecture notes by R. Askey \cite{A}. Writing $B(x,y)$ for the Beta function, they showed that
\begin{align}\label{8pi3}
    \int_{-\infty}^{\infty} t \sinh \pi t  \prod_{\epsilon_1,\epsilon_2=\pm 1} B(\tfrac{1}{4}+\epsilon_1 i\tfrac{t}{2}+ \epsilon_2 iT, \tfrac{1}{4}- \epsilon_1 i\tfrac{t}{2}) dt = 8 \pi^3.
\end{align}
With Kuznetsov's formula (see \cite[Thm 9.3]{I2}), this identity immediately yields the following lemma.

\begin{lemma}\label{kuz}
We have
\begin{align*}
    \sum_{u\in \mathcal{B}(N)}|\rho_u|^2 \lambda_u(m)\lambda_u(n)  \prod_{\epsilon_1,\epsilon_2=\pm 1} B(\tfrac{1}{4}+\epsilon_1 i\tfrac{t_u}{2}+ \epsilon_2 iT, \tfrac{1}{4}- \epsilon_1 i\tfrac{t_u}{2}) = \delta_{m=n} 8\pi + \text{Off-Diagonal Terms}.
\end{align*}
\end{lemma}
\begin{remark}
The off-diagonal terms can be ignored for the main term of $I_1$, due to the Recipe.
\end{remark}

\section{Proof of Theorem \ref{main}}
\subsection{Initial cleanings}
Recall notations in Section \ref{orthobase}. Let $E=E_{\infty}(z,\frac{1}{2}+iT,\chi)$ be a newform Eisenstein series of even primitive nebentypus mod $N$. We begin with simplifying
\begin{align}\label{withEis}
    \langle |E|^2 - \mathcal{E}, |E|^2 - \mathcal{E} \rangle = \sum_{u\in \mathcal{O}(N)} \langle |E|^2,u\rangle \langle u,|E|^2\rangle.
\end{align}
When $N$ is prime, \eqref{BM} reduces the above sum to
\begin{align*}
    \sum_{u\in \mathcal{B}(N)} \langle |E|^2,u\rangle \langle u,|E|^2\rangle + \sum_{u\in \mathcal{B}(1)} \sum_{\ell \mid N} \langle |E|^2,u^{<\ell>}\rangle \langle u^{<\ell>},|E|^2\rangle.
\end{align*}
Then by \eqref{oldclass} for all $u\in \mathcal{O}(1)$ we have
\begin{align*}
    \langle |E|^2,u^{<1>}\rangle &= \langle |E|^2,u\rangle\\
    \langle |E|^2,u^{<N>}\rangle &= \langle |E|^2, \xi_N(1) u + \xi_N(N) u|_N \rangle = \overline{(\xi_N(1) + \xi_N(N))} \langle |E|^2,u\rangle,
\end{align*}
so \eqref{BMsharp} implies
\begin{align*}
    \sum_{u\in \mathcal{B}(1)} \sum_{\ell \mid N} \langle |E|^2,u^{<\ell>}\rangle \langle u^{<\ell>},|E|^2\rangle = (2 + O(N^{-1+2\theta+\varepsilon})) \sum_{u\in \mathcal{B}(1)} \langle |E|^2,u\rangle \langle u,|E|^2\rangle.
\end{align*}
With \eqref{expression}, the spectral sum over $u\in \mathcal{B}(1)$ can be bounded by $N^{-2+\varepsilon}$ by taking the Lindel\"of bound for the L-functions, which is negligible under comparison with the contribution from $\mathcal{B}(N)$. Thus
\begin{align*}
    \langle |E|^2 - \mathcal{E}, |E|^2 - \mathcal{E} \rangle \sim \sum_{u\in \mathcal{B}(N)} \langle |E|^2,u\rangle \langle u,|E|^2\rangle.
\end{align*}
\begin{remark}\label{relaxation}
Following similar argument we can see the general case. For $N>1$, we have
\begin{align*}
    \sum_{u\in \mathcal{O}(N)} \langle |E|^2,u\rangle \langle u,|E|^2\rangle = \sum_{N=ML} \sum_{u\in \mathcal{B}(M)} \sum_{\ell \mid L} \langle |E|^2,u^{<\ell>}\rangle \langle u^{<\ell>},|E|^2\rangle.
\end{align*}
Here $u^{<\ell>}$ can similarly be expressed as a linear combination of $u|_d$ for $d\mid \ell$, and the sum over all linear coefficients is $O(\ell^{\epsilon})$, due to \cite{BM}.
For all $M<N$ and $L>1$, we have
\begin{align*}
    \sum_{u\in \mathcal{B}(M)} \langle |E|^2,u\rangle \langle u,|E|^2\rangle \ll \frac{N^{\varepsilon}}{\nu(N)} \frac{M}{N}.
\end{align*}
In total there is
\begin{align*}
    \sum_{u\in \mathcal{O}(N)} \langle |E|^2,u\rangle \langle u,|E|^2\rangle = \frac{N^{\varepsilon}}{\nu(N)} \sum_{\substack{M\mid N \\ M<N}} \frac{M}{N} + \sum_{u\in \mathcal{B}(N)} \langle |E|^2,u\rangle \langle u,|E|^2\rangle.
\end{align*}
As the proof goes, the second term is of size $\nu^{-1}(N) \log^2 N$, so we can allow $N$ to satisfy 
\begin{align*}
    \sum_{\substack{M\mid N \\ M>1}} \frac{1}{M} \ll N^{-\delta},
\end{align*}
for any fixed $\delta>0$, as is claimed in Remark \ref{composite}.
\end{remark}

\subsection{Product of formal approximate functional equations}
For convenience, denote
\begin{align*}
\chi^{+1}=\chi, \quad \chi^{-1}=\overline{\chi}.
\end{align*}
Also write $\epsilon=(\epsilon_1,\epsilon_2,\epsilon_3,\epsilon_4) \in \{\pm 1\}^4$, $\alpha=(\alpha_1,\alpha_2,\alpha_3,\alpha_4) \in \mathbb{C}^4$ and $\alpha_0=(0,0,2iT,-2iT)$.
Insert the formal approximate functional equations to \eqref{expression} and sum over $u\in \mathcal{B}(N)$. The recipe says
\begin{multline} \label{limit}
    \sum_{u} |\langle |E|^2, u \rangle|^2 \sim \frac{1}{8N} \lim_{\alpha\rightarrow\alpha_0} \sum_{u} (1+\lambda_u(-1))  W_u \sum_{\epsilon_1,...,\epsilon_4 = \pm 1} \prod_{j=1}^4 X_j f_j
    \\
    \sum_{n_1,...,n_4 \geq 1} \frac{\lambda_u(n_1)\lambda_u(n_2)\lambda_u(n_3)\lambda_u(n_4)\chi^{\epsilon_3}(n_3)\overline{\chi}^{\epsilon_4}(n_4)}{n_1^{\frac{1}{2}+\epsilon_1 \alpha_1} n_2^{\frac{1}{2}+\epsilon_2\alpha_2} n_3^{\frac{1}{2}+\epsilon_3\alpha_3} n_4^{\frac{1}{2}+\epsilon_4\alpha_4}},
\end{multline}
with
\begin{align}\label{Wu}
    W_u = \frac{|\rho_u|^2 }{ |L(1+2iT,\chi)|^4} \prod_{\epsilon_1,\epsilon_2=\pm 1} B(\tfrac{1}{4}+\epsilon_1 i\tfrac{t}{2}+ \epsilon_2 iT, \tfrac{1}{4}- \epsilon_1 i\tfrac{t}{2}),
\end{align}
$X_j=X_j(\epsilon_j)$ and $f_j=f_j(\epsilon_j)$ for $j=1,2,3,4$. Moreover,
\begin{align*}
    X_j(\epsilon_j) = \Big(\lambda(-1&)\Big)^{\frac{1-\epsilon_j}{2}}, j=1,2; \quad 
    X_3(\epsilon_3) = \Big(\lambda_u(-1) (\frac{\tau(\chi)}{\sqrt{N}})^2\Big)^{\frac{1-\epsilon_3}{2}}, \quad X_4(\epsilon_4) =\Big(\lambda_u(-1) (\frac{\tau(\overline{\chi})}{\sqrt{N}})^2\Big)^{\frac{1-\epsilon_4}{2}} \\
    f_j(\epsilon_j) &= \Big(f(\frac{1}{2}+\alpha_j,u)\Big)^{\frac{1-\epsilon_j}{2}}, j=1,2; \quad
    f_j(\epsilon_j) = \Big(f(\tfrac{1}{2}+\alpha_j,u\otimes\chi)\Big)^{\frac{1-\epsilon_j}{2}}, j=3,4.
\end{align*}

\begin{remark}\label{anypath}
Since the right hand side in \eqref{limit} is analytic around $\alpha=\alpha_0$, we can take any path to approach this point in a neighborhood of it. This is useful in the last step of the proof.
\end{remark}

By Hecke relation, the multi-Dirichlet series equals
\begin{align*}
    \zeta(1+\epsilon_1\alpha_1+\epsilon_2\alpha_2) L(1+\epsilon_3\alpha_3+\epsilon_4\alpha_4, \chi_{0,N}) \sum_{n_1,...,n_4 \geq 1}\frac{\lambda_u(n_1 n_2) \lambda_u(n_3 n_4) \chi^{\epsilon_3}(n_3)\overline{\chi}^{\epsilon_4}(n_4)}{n_1^{\frac{1}{2}+\epsilon_1 \alpha_1} n_2^{\frac{1}{2}+\epsilon_2\alpha_2} n_3^{\frac{1}{2}+\epsilon_3\alpha_3} n_4^{\frac{1}{2}+\epsilon_4\alpha_4}},
\end{align*}

With \eqref{withEis} and \eqref{limit}, we can write
\begin{align}\label{milestone1}
    \langle |E|^2 - \mathcal{E}, |E|^2 - \mathcal{E} \rangle = \lim_{\alpha\rightarrow\alpha_0}  \sum_{\epsilon_1,...,\epsilon_4=\pm 1} S(\epsilon,\alpha),
\end{align}
where $S(\epsilon,\alpha)$ equals
\begin{align}\label{ABC}
    \sum_{u} \underbrace{(1+\lambda_u(-1)) \prod_{j=1}^4 X_j}_{A_u} \underbrace{W_u \prod_{j=1}^4 f_j  \sum_{n_1,...,n_4 \geq 1} \frac{\lambda_u(n_1)\lambda_u(n_2)\lambda_u(n_3)\lambda_u(n_4)\chi^{\epsilon_3}(n_3)\overline{\chi}^{\epsilon_4}(n_4)}{n_1^{\frac{1}{2}+\epsilon_1 \alpha_1} n_2^{\frac{1}{2}+\epsilon_2\alpha_2} n_3^{\frac{1}{2}+\epsilon_3\alpha_3} n_4^{\frac{1}{2}+\epsilon_4\alpha_4}} }_{B_u}.
\end{align}
Since
\begin{align*}
    \frac{1}{|\mathcal{B}(N)|} \sum_{u \in \mathcal{B}(N)} (1+\lambda_u(-1)) \prod_{j=1}^4 X_j = (\frac{\tau(\chi)}{\sqrt{N}})^{\epsilon_3-\epsilon_4} + o(1),
\end{align*}
following the recipe, we replace $A_u$ with $(\tfrac{\tau(\chi)}{\sqrt{N}})^{\epsilon_3-\epsilon_4}$, its ``expected mean'' in the family, in \eqref{ABC}.

\subsection{Reduction of the multi-Dirichlet series}
Now we see that $S(\epsilon,\alpha)$ asymptotically equals
\begin{align*}
    \frac{1}{8N}  (\frac{\tau(\chi)}{\sqrt{N}})^{\epsilon_3-\epsilon_4} \sum_{u}  W_u \prod_{j=1}^4 f_j
    \sum_{n_1,...,n_4 \geq 1} \frac{\lambda_u(n_1)\lambda_u(n_2)\lambda_u(n_3)\lambda_u(n_4)\chi^{\epsilon_3}(n_3)\overline{\chi}^{\epsilon_4}(n_4)}{n_1^{\frac{1}{2}+\epsilon_1 \alpha_1} n_2^{\frac{1}{2}+\epsilon_2\alpha_2} n_3^{\frac{1}{2}+\epsilon_3\alpha_3} n_4^{\frac{1}{2}+\epsilon_4\alpha_4}}.
\end{align*}
Applying Hecke relation, we can rewrite the quadruple sum on the right hand side as
\begin{align*}
    \zeta(1+\epsilon_1\alpha_1+\epsilon_2\alpha_2) L(1+\epsilon_3\alpha_3+\epsilon_4\alpha_4, \chi^{\epsilon_3}\overline{\chi}^{\epsilon_4}) \sum_{n_1,...,n_4 \geq 1}\frac{\lambda_u(n_1 n_2) \lambda_u(n_3 n_4) \chi^{\epsilon_3}(n_3)\overline{\chi}^{\epsilon_4}(n_4)}{n_1^{\frac{1}{2}+\epsilon_1 \alpha_1} n_2^{\frac{1}{2}+\epsilon_2\alpha_2} n_3^{\frac{1}{2}+\epsilon_3\alpha_3} n_4^{\frac{1}{2}+\epsilon_4\alpha_4}}.
\end{align*}
For convenience, denote
\begin{align}\label{hu}
    h_u(\epsilon,\alpha,t) = \prod_{j=1}^4 f_j  \Big/ \Big(  (\frac{\pi^2}{N})^{\frac{1-\epsilon_1}{2}\alpha_1+\frac{1-\epsilon_2}{2}\alpha_2} (\frac{\pi^2}{N^2})^{\frac{1-\epsilon_3}{2}\alpha_3+\frac{1-\epsilon_4}{2}\alpha_4} \Big)
\end{align}
for the product of Gamma fractions. Note $f_j$ depends on $t=t_u$.
Then
\begin{multline}\label{ss}
    S(\epsilon,\alpha) \sim \frac{1}{8N} (\frac{\tau(\chi)}{\sqrt{N}})^{\epsilon_3-\epsilon_4}  (\frac{\pi^2}{N})^{\frac{1-\epsilon_1}{2}\alpha_1+\frac{1-\epsilon_2}{2}\alpha_2} (\frac{\pi^2}{N^2})^{\frac{1-\epsilon_3}{2}\alpha_3+\frac{1-\epsilon_4}{2}\alpha_4} \zeta(1+\epsilon_1\alpha_1+\epsilon_2\alpha_2) \\
     L(1+\epsilon_3\alpha_3+\epsilon_4\alpha_4, \chi^{\epsilon_3}\overline{\chi}^{\epsilon_4})
    \sum_{u} h_u W_u \sum_{n_1,...,n_4 \geq 1}\frac{\lambda_u(n_1 n_2) \lambda_u(n_3 n_4) \chi^{\epsilon_3}(n_3)\overline{\chi}^{\epsilon_4}(n_4)}{n_1^{\frac{1}{2}+\epsilon_1 \alpha_1} n_2^{\frac{1}{2}+\epsilon_2\alpha_2} n_3^{\frac{1}{2}+\epsilon_3\alpha_3} n_4^{\frac{1}{2}+\epsilon_4\alpha_4}}.
\end{multline}
The next step of the recipe, we need to find the diagonal term for
\begin{align*}
    \sum_{u} h_u W_u \sum_{n_1,...,n_4 \geq 1}\frac{\lambda_u(n_1 n_2) \lambda_u(n_3 n_4) \chi^{\epsilon_3}(n_3)\overline{\chi}^{\epsilon_4}(n_4)}{n_1^{\frac{1}{2}+\epsilon_1 \alpha_1} n_2^{\frac{1}{2}+\epsilon_2\alpha_2} n_3^{\frac{1}{2}+\epsilon_3\alpha_3} n_4^{\frac{1}{2}+\epsilon_4\alpha_4}}
\end{align*}
in the Kuznetsov trace formula. In general,
\begin{align*}
    \sum_{u} h_u W_u \lambda_u(n_1 n_2) \lambda_u(n_3 n_4) \sim \delta_{n_1n_1=n_3n_4} F_{\epsilon}(h_u;T) |L(1+2iT,\chi)|^{-4},
\end{align*}
with
\begin{align*}
    F_{\epsilon}(h_u;T) =
    \frac{1}{\pi^2} \int_{-\infty}^{\infty} t \sinh(\pi t)
    \prod_{\epsilon_1,\epsilon_2=\pm 1} B(\tfrac{1}{4}+\epsilon_1 i\tfrac{t_u}{2}+ \epsilon_2 iT, \tfrac{1}{4}- \epsilon_1 i\tfrac{t_u}{2})
    h_u(\epsilon,\alpha, t) dt,
\end{align*}
and by Lemma \ref{kuz}, $F_{\epsilon}(h_u;T) =8\pi$ when $\epsilon_3=\epsilon_4$. If $\epsilon_3 \neq \epsilon_4$,
as $T\rightarrow 0$ we also have $F_{\epsilon}(h_u;T)\rightarrow 8\pi$.

For the remaining quadruple sum
\begin{align*}
    \sum_{\substack{n_1,...,n_4 \geq 1 \\ n_1n_2=n_3n_4}}\frac{ \chi^{\epsilon_3}(n_3)\overline{\chi}^{\epsilon_4}(n_4)}{n_1^{\frac{1}{2}+\epsilon_1 \alpha_1} n_2^{\frac{1}{2}+\epsilon_2\alpha_2} n_3^{\frac{1}{2}+\epsilon_3\alpha_3} n_4^{\frac{1}{2}+\epsilon_4\alpha_4}},
\end{align*}
we can apply Ramanujan's formula for twisted Dirichlet series (see \cite[(13.1)]{I1}) to rewrite it as
\begin{align*}
    \frac{L(1+\epsilon_1\alpha_1+\epsilon_3\alpha_3,\chi^{\epsilon_3}) L(1+\epsilon_2\alpha_2+\epsilon_3\alpha_3, \chi^{\epsilon_3}) L(1+\epsilon_1\alpha+\epsilon_4\alpha_4, \overline{\chi}^{\epsilon_4}) L(1+\epsilon_2\alpha_2+\epsilon_4\alpha_4, \overline{\chi}^{\epsilon_4})}{L(2+\epsilon_1\alpha_1+\epsilon_2\alpha_2+\epsilon_3\alpha_3+\epsilon_4\alpha_4, \chi^{\epsilon_3}\overline{\chi}^{\epsilon_4})}.
\end{align*}
To summarize, \eqref{ss} can be rewritten as
\begin{multline}\label{milestone2}
    S(\epsilon,\alpha) \sim \frac{F_{\epsilon}(h;T)}{8N |L(1+2iT,\chi)|^4} (\frac{\tau(\chi)}{\sqrt{N}})^{\epsilon_3-\epsilon_4}
    \prod_{j=1}^4 \pi^{(1-\epsilon_j)\alpha_j} \\
    N^{\frac{\epsilon_1 -1}{2}\alpha_1+\frac{\epsilon_2 -1}{2}\alpha_2+(\epsilon_3 -1)\alpha_3 +(\epsilon_4 -1)\alpha_4}
    \zeta(1+\epsilon_1\alpha_1+\epsilon_2\alpha_2)
     L(1+\epsilon_3\alpha_3+\epsilon_4\alpha_4, \chi^{\epsilon_3}\overline{\chi}^{\epsilon_4}) \\
     \frac{\prod_{j=1,2}L(1+\epsilon_j\alpha_j+\epsilon_3\alpha_3,\chi^{\epsilon_3})   L(1+\epsilon_j\alpha_j+\epsilon_4\alpha_4, \overline{\chi}^{\epsilon_4})}{L(2+\epsilon_1\alpha_1+\epsilon_2\alpha_2+\epsilon_3\alpha_3+\epsilon_4\alpha_4, \chi^{\epsilon_3}\overline{\chi}^{\epsilon_4})}.
\end{multline}

\subsection{Final evaluation}
Recall Remark \ref{anypath}, to calculate $\lim_{\alpha\rightarrow\alpha_0}\sum_\epsilon S(\epsilon,\alpha)$, we can follow such an order:
\begin{itemize}
    \item Substitute $\alpha=(0,\eta',2iT,-2iT+\eta)$ for $\eta>\eta'>0$.
    \item Compute $\lim_{\eta'\rightarrow 0} \sum_{\epsilon_1,\epsilon_2=\pm 1}S(\epsilon,\alpha)$. Write the result as $R_{\epsilon_3,\epsilon_4}(\eta)$.
    \item Compute $\lim_{\eta \rightarrow 0} \sum_{\epsilon_3,\epsilon_4=\pm 1}R_{\epsilon_3,\epsilon_4}(\eta)$.
\end{itemize}
Specifically, as $\eta'\rightarrow 0$, we have
\begin{multline*}
     F_{\epsilon}(h;T) \pi^{(1-\epsilon_2)\eta'}N^{\frac{\epsilon_2 -1}{2}\eta'} \zeta(1+\epsilon_2\eta')
     \frac{L(1+\epsilon_2\eta'+\epsilon_3 2iT, \chi^{\epsilon_3}) L(1+\epsilon_2\eta'+\epsilon_4(-2iT+\eta), \overline{\chi}^{\epsilon_4})}{L(2+\epsilon_2\eta'+(\epsilon_3-\epsilon_4)2iT+\epsilon_4\eta, \chi^{\epsilon_3}\overline{\chi}^{\epsilon_4})} \\
    =\frac{K_{\epsilon_3,\epsilon_4}(T,\eta)}{\epsilon_2\eta'} - \delta_{\epsilon_2=-1} H_{\epsilon}(\eta) \log N 
    \frac{L(1+\epsilon_3 2iT, \chi^{\epsilon_3}) L(1+\epsilon_4(-2iT+\eta), \overline{\chi}^{\epsilon_4})}{L(2+(\epsilon_3-\epsilon_4)2iT+\epsilon_4\eta, \chi^{\epsilon_3}\overline{\chi}^{\epsilon_4})} + O(\log\log N) + O(\eta'),
\end{multline*}
where $H_{\epsilon}(\eta)=F_{\epsilon}(h_0;T)$ with $h_0=h_u(\epsilon,(0,0,2iT,-2iT+\eta),t)$, and
\begin{align*}
    K_{\epsilon_3,\epsilon_4}(T,\eta) = \frac{L(1+\epsilon_3 2iT, \chi^{\epsilon_3}) L(1+\epsilon_4(-2iT+\eta), \overline{\chi}^{\epsilon_4})}{L(2+(\epsilon_3-\epsilon_4)2iT+\epsilon_4\eta, \chi^{\epsilon_3}\overline{\chi}^{\epsilon_4})}.
\end{align*}
Since $K_{\epsilon_3,\epsilon_4}(T,\eta)$ is independent of $\epsilon_1,\epsilon_2$, and so is $H_{\epsilon}(\eta)$ when $\alpha_1=0, \alpha_2\rightarrow 0$, $R_{\epsilon_3,\epsilon_4}(\eta)$ equals
\begin{multline*}
     \lim_{\eta'\rightarrow 0}\sum_{\epsilon_1,\epsilon_2 =\pm 1} S(\epsilon,\alpha) = \frac{\pi^{(\epsilon_4-\epsilon_3)2iT }N^{(\epsilon_3 -\epsilon_4)2iT}}{8N |L(1+2iT,\chi)|^4}
     (\frac{\tau(\chi)}{\sqrt{N}})^{\epsilon_3-\epsilon_4}
    L(1+(\epsilon_3-\epsilon_4)2iT+\epsilon_4\eta, \chi^{\epsilon_3}\overline{\chi}^{\epsilon_4}) 
    \\
    \pi^{(1-\epsilon_4)\eta} N^{(\epsilon_4-1)\eta}
    \frac{L^2(1+\epsilon_3 2iT,\chi^{\epsilon_3}) L^2(1+\epsilon_4 (-2iT+\eta),\overline{\chi}^{\epsilon_4})}{L(2+(\epsilon_3-\epsilon_4)2iT+\epsilon_4\eta, \chi^{\epsilon_3}\overline{\chi}^{\epsilon_4})}
    \Big( -2 H_{\epsilon}(\eta) \log N +O(\log\log N) \Big).
\end{multline*}

Now we analyze the Laurent expansion of $R_{\epsilon_3,\epsilon_4}(\eta)$ around $\eta=0$. We need to discuss $4$ cases:

\noindent \textbf{Case 1.} If $\epsilon_3\neq\epsilon_4$ and $\chi$ is complex, then $L(1+(\epsilon_3-\epsilon_4)2iT+\epsilon_4\eta, \chi^{\epsilon_3}\overline{\chi}^{\epsilon_4})$ is not principal, so $R_{\epsilon_3,\epsilon_4}(\eta)$ is analytic at $\eta=0$, and by GRH, for any fixed $T\in \mathbb{R}$ we have (note $H_{\epsilon}(0)=O(1)$)
\begin{align}\label{nopole}
    R_{\epsilon_3,\epsilon_4}(0) = O(\log N (\log\log N)^9). 
\end{align}
\noindent \textbf{Case 2.} If $\epsilon_3\neq\epsilon_4$, $\chi$ is quadratic, but $T\neq 0$, then there is no pole around $\eta=0$ for the principal L-function $L(1+(\epsilon_3-\epsilon_4)2iT+\epsilon_4\eta, \chi^{\epsilon_3}\overline{\chi}^{\epsilon_4})$, so we similarly have \eqref{nopole}.

\begin{remark}\label{varyTA}
As we can see in discussion for Case 4, the value of $T$ no longer matters when $\epsilon_3=\epsilon_4$, so we explain our claim in Remark \ref{varyTQ} here. When $\chi$ is quadratic and $T=T(N)$ tends to $0$, the principal L-function $L(1+(\epsilon_3-\epsilon_4)2iT+\epsilon_4\eta, \chi^{\epsilon_3}\overline{\chi}^{\epsilon_4})$ equals
\begin{align*}
    \zeta(1+(\epsilon_3-\epsilon_4) 2iT +\epsilon_4 \eta) \prod_{p\mid N} (1-\frac{1}{p^{1+(\epsilon_3-\epsilon_4) 2iT +\epsilon_4 \eta}})
\end{align*}
and has large values with small $T$ and $\eta$. Therefore, we have (recall $H_{\epsilon}(0)=8\pi$ for all $\epsilon$)
\begin{align*}
    R_{+1,-1}(0) = \frac{\pi^{-4iT}N^{4iT}}{8N \zeta(2+4iT)}
    \zeta(1+4iT) \prod_{p\mid N} \frac{1-p^{-(1+4iT)}}{1-p^{-(2+4iT)}}
    \Big( -16 \pi \log N +O(\log\log N) \Big)
\end{align*}
and
\begin{align*}
    R_{-1,+1}(0) = \frac{\pi^{4iT}N^{-4iT}}{8N \zeta(2+4iT)}
    \zeta(1-4iT) \prod_{p\mid N} \frac{1-p^{-(1-4iT)}}{1-p^{-(2-4iT)}}
    \Big( -16 \pi\log N +O(\log\log N) \Big).
\end{align*}
Consequently,
\begin{align}\label{upndown}
    R_{+1,-1}(0) + R_{-1,+1}(0) = \frac{-\log N}{4\nu(N) \zeta(2)} \Big( \frac{e^{4iT \log N}}{4iT} + \frac{e^{-4iT \log N}}{-4iT} \Big) + O(\log N (\log\log N)^3).
\end{align}
If $4iT \log N = o(1)$, then $2\log N$ approximates the factor in the parentheses with an error up to $O(|T| \log^2 N)$. So, for all $|T|\ll \log^{-1-\delta} N$, we have $R_{+1,-1}(0) + R_{-1,+1}(0)=\frac{24}{\pi} \nu^{-1}(N) \log^2 N$.
On the other hand, by \eqref{upndown}, $R_{+1,-1}(0) + R_{-1,+1}(0)$ has trivial bound $\nu^{-1}(N)|T|^{-1} \log N$. So, for whatever $|T| \gg \log^{-1+\delta}N$, $R_{+1,-1}(0) + R_{-1,+1}(0)$ will be considerably less than $\nu^{-1}(N) \log^2 N$.
\end{remark} 

\noindent \textbf{Case 3.} If $\epsilon_3\neq\epsilon_4$, $\chi$ is quadratic, and $T= 0$, then $L(1+(\epsilon_3-\epsilon_4)2iT+\epsilon_4\eta, \chi^{\epsilon_3}\overline{\chi}^{\epsilon_4})$ has a pole at $\eta=0$. So $R_{\epsilon_3,\epsilon_4}(\eta)$ equals
\begin{align*}
     R_{\epsilon_3,\epsilon_4}(\eta) = \frac{\pi^{(1-\epsilon_4)\eta} N^{(\epsilon_4-1)\eta}}{8N} H_{\epsilon}(\eta) \prod_{p\mid N}\frac{1-p^{-1-\epsilon_4\eta}}{1-p^{-2-\epsilon_4\eta}}\frac{\zeta(1+\epsilon_4 \eta)}
    {\zeta(2+\epsilon_4\eta)} \Big( -2 H_{\epsilon}(\eta) \log N +O(\log\log N) \Big), 
\end{align*}
and the Laurent expansion around $\eta=0$ is (note $\zeta(2)=\tfrac{\pi^2}{6}$ and $\nu(N)=N\prod_{p\mid N}(1+p^{-1})$)
\begin{multline*}
    \frac{6}{\nu(N) \pi}
    \Big( \frac{1}{\epsilon_4 \eta} + (\epsilon_4-1) \log N +O\big((\log\log N)^5 \big) +O(\eta) \Big)
    \Big( -2 \log N +O(\log\log N) \Big) \\
    = \delta_{\epsilon_4=-1} \frac{24}{\nu(N) \pi} \log^2 N +O(\log N (\log\log N)^5).
\end{multline*}
\noindent \textbf{Case 4.} If $\epsilon_3=\epsilon_4$, we have 
\begin{multline}\label{milestone3}
    R_{\epsilon_3,\epsilon_4}(\eta) = \frac{\pi^{(1-\epsilon_4)\eta} N^{(\epsilon_4-1)\eta}}{8N} \frac{L^2(1+\epsilon_3 2iT,\chi^{\epsilon_3}) L^2(1+\epsilon_4 (-2iT+\eta),\overline{\chi}^{\epsilon_4})}{|L(1+2iT,\chi)|^4} H_{\epsilon}(\eta)   \\
    \prod_{p\mid N}\frac{1-p^{-1-\epsilon_4\eta}}{1-p^{-2-\epsilon_4\eta}}
    \frac{\zeta(1+\epsilon_4 \eta)}
    {\zeta(2+\epsilon_4\eta)} \Big( -2 H_{\epsilon}(\eta) \log N +O(\log\log N) \Big),
\end{multline}
which further equals (note here $H_{\epsilon}(0)=8\pi$ and $H_{\epsilon}'(0)=O(1)$)
\begin{multline*}
    \frac{6}{\nu(N) \pi} \frac{L^2(1+\epsilon_3 2iT,\chi^{\epsilon_3}) L^2(1+\epsilon_4 (-2iT),\overline{\chi}^{\epsilon_4})}{|L(1+2iT,\chi)|^4} \\
    \Big( \frac{1}{\epsilon_4 \eta} + (\epsilon_4-1) \log N +O\big((\log\log N)^5 \big) +O(\eta) \Big)
    \Big( -2 \log N +O(\log\log N) \Big).
\end{multline*}
With above calculations we can rewrite it as
\begin{align*}
 \delta_{\epsilon_4=-1} \frac{24}{\nu(N)\pi} \log^2 N + O(\log N (\log\log N)^9).
\end{align*}

\noindent Summing up, we have
\begin{align*}
    I_1 \sim \lim_{\eta\rightarrow 0}  \sum_{\epsilon_3,\epsilon_4 =\pm 1} R_{\epsilon_3,\epsilon_4}(\eta) = \frac{24}{\nu(N)\pi}
    \Big( 1 + \delta_{T=0}\delta_{\chi \text{ quadratic}} \Big) \log^2 N.
\end{align*}

\section{Proof of Theorem \ref{0diff}}
The proof follows the lines of Djankovi\'c and Khan in proving the $t$-aspect formula
\begin{align*}
    \langle |E_t|^2, |E_t|^2 \rangle_{\mathrm{reg}} - \langle |\Lambda^Y E_t|^2, |\Lambda^Y E_t| \rangle \sim \frac{36}{\pi} \log^2 t,
\end{align*}
provided that $1<Y\ll \log t$. Here $E_t=E(z,\tfrac{1}{2}+it)$ is the classical Eisenstein series. 

For Atkin-Lehner cusp $\mathfrak{a}$, Define $e_{\mathfrak{a}}$ to be the main term of $E_{\mathfrak{a}}$ as $z$ approaches each cusp. Alternatively, we have for $s=\tfrac{1}{2}$ and $\chi^2=\chi_{0,N}$
\begin{align}\label{ecusp}
    e_{\mathfrak{a}}|_{\sigma_{\mathfrak{b}}} =
    \begin{cases}
    \sqrt{y} & \text{ if } \mathfrak{b=a,a^*} \\
    0 & \text{ otherwise}.
    \end{cases}
\end{align}

\noindent Following the definitions, we have
\begin{multline}\label{broad}
    \langle |E|^2, |E|^2 \rangle_{\mathrm{reg}} - \langle |\Lambda^Y E|^2, |\Lambda^Y E|^2 \rangle = \int_{\mathcal{F}} e_{\mathfrak{a}}^2 \overline{E^Y_{\mathfrak{a}}}^2 + \overline{e_{\mathfrak{a}}}^2 (E^Y_{\mathfrak{a}})^2 + 4|e_{\mathfrak{a}}|^2 |E^Y_{\mathfrak{a}}|^2 d\mu \\
    2\int_{\mathcal{F}} \Big( |E^Y_{\mathfrak{a}}|^2 \overline{E^Y_{\mathfrak{a}}} e_{\mathfrak{a}} + |E^Y_{\mathfrak{a}}|^2 E^Y_{\mathfrak{a}} \overline{e_{\mathfrak{a}}} \Big) d\mu + \Phi(Y),
\end{multline}
where $\mathcal{F}=\cup_{\mathfrak{a}} \sigma_{\mathfrak{a}} \mathcal{F}_{\infty}(Y)$ is the cuspical zone of height $Y$, with scaling matrix $\sigma_{\mathfrak{a}}$, $\Phi(Y)$ is some $N$-independent function in $Y$, and $\mathcal{F}_{\infty}(Y)=\{(x,y) | 0<1\leq 1, y>Y \}$.

\subsection{Initial cleanings}
For primitive $\chi_1$, $\chi_2$, write $\psi=\chi_1\chi_2$. It is easy to see if $\chi=\chi\overline{\chi_2}$ is even, then so is $\psi$. Denote the completed Eisenstein series $E_{\chi_1,\chi_2}^*(z,s)=\theta_{\chi_1,\chi_2}(s) E_{\chi_1,\chi_2}(z,s)$, with $\theta_{\chi_1,\chi_2}(\tfrac{1}{2})=\Lambda^{-1}(1,\psi)$ as in \eqref{theta}.
Write out the completion factor with Young's formula \cite[(9.1)]{Young}:
\begin{align*}
    E_{\mathfrak{a}}(z,\tfrac{1}{2},\chi)= N^{-\frac{1}{2}} E_{\chi_1,\chi_2}(z,\tfrac{1}{2}) = N^{-\frac{1}{2}} \theta^{-1}_{\chi_1,\chi_2}(\tfrac{1}{2}) E^*_{\chi_1,\chi_2}(z,\tfrac{1}{2}),
\end{align*}
and
\begin{align}\label{Eslash}
    E_{\mathfrak{a}}(\sigma_{\mathfrak{a}} z, \tfrac{1}{2},\chi) = E_{\mathfrak{a}^*}(\sigma_{\mathfrak{a}} z, \tfrac{1}{2},\chi) = \pm N^{-\frac{1}{2}} \theta^{-1}_{\chi_1,\chi_2}(\tfrac{1}{2}) E^*_{1,\psi}(z,\tfrac{1}{2}).
\end{align}

\noindent When $s=\tfrac{1}{2}$ and $\chi$ is quadratic, $E_{\mathfrak{a}}(z,s,\chi)$ is real-valued, and so is its $Y$-truncation. Then by \eqref{ecusp}

\begin{multline*}
    \int_{\mathcal{F}} e_{\mathfrak{a}}^2 \overline{E^Y_{\mathfrak{a}}}^2 + \overline{e_{\mathfrak{a}}}^2 (E^Y_{\mathfrak{a}})^2 + 4|e_{\mathfrak{a}}|^2 |E^Y_{\mathfrak{a}}|^2 d\mu = \sum_{\mathfrak{b}} \int_Y^{\infty}\int_0^1 \Big( e_{\mathfrak{a}}^2 \overline{E^Y_{\mathfrak{a}}}^2 + \overline{e_{\mathfrak{a}}}^2 (E^Y_{\mathfrak{a}})^2 + 4|e_{\mathfrak{a}}|^2 |E^Y_{\mathfrak{a}}|^2 \Big)\Big|_{\sigma_{\mathfrak{b}}^{-1}} d\mu \\
    =\int_Y^{\infty}\int_0^1 \Big( e_{\mathfrak{a}}^2 \overline{E^Y_{\mathfrak{a}}}^2 + \overline{e_{\mathfrak{a}}}^2 (E^Y_{\mathfrak{a}})^2 + 4|e_{\mathfrak{a}}|^2 |E^Y_{\mathfrak{a}}|^2 \Big)\Big|_{\sigma_{\mathfrak{a}}^{-1}} + \Big( e_{\mathfrak{a}}^2 \overline{E^Y_{\mathfrak{a}}}^2 + \overline{e_{\mathfrak{a}}}^2 (E^Y_{\mathfrak{a}})^2 + 4|e_{\mathfrak{a}}|^2 |E^Y_{\mathfrak{a}}|^2 \Big)\Big|_{\sigma_{\mathfrak{a}^*}^{-1}} d\mu.
\end{multline*}
Since $\mathfrak{a}$ and $\mathfrak{a}^*$ are Atkin-Lehner, their scaling matrices are some matrix involution. Thus we can write the integral as
\begin{align}\label{DK5.1}
    \int_Y^{\infty}\int_0^1 \Big( e_{\mathfrak{a}}^2 \overline{E^Y_{\mathfrak{a}}}^2 + \overline{e_{\mathfrak{a}}}^2 (E^Y_{\mathfrak{a}})^2 + 4|e_{\mathfrak{a}}|^2 |E^Y_{\mathfrak{a}}|^2 \Big)\Big|_{\sigma_{\mathfrak{a}}} + \Big( e_{\mathfrak{a}}^2 \overline{E^Y_{\mathfrak{a}}}^2 + \overline{e_{\mathfrak{a}}}^2 (E^Y_{\mathfrak{a}})^2 + 4|e_{\mathfrak{a}}|^2 |E^Y_{\mathfrak{a}}|^2 \Big)\Big|_{\sigma_{\mathfrak{a}^*}} d\mu.
\end{align}

\noindent On the other hand, by \eqref{ecusp} and \eqref{Eslash} we have
\begin{align*}
    (e_{\mathfrak{a}}|_{\sigma_{\mathfrak{a}}})^2 = (e_{\mathfrak{a}}|_{\sigma_{\mathfrak{a}^*}})^2 = (\overline{e_{\mathfrak{a}}}|_{\sigma_{\mathfrak{a}}})^2= (\overline{e_{\mathfrak{a}}}|_{\sigma_{\mathfrak{a}^*}})^2 = |e_{\mathfrak{a}}|_{\sigma_{\mathfrak{a}}}|^2 = |e_{\mathfrak{a}}|_{\sigma_{\mathfrak{a}^*}}|^2 &=y, \\
    (E^Y_{\mathfrak{a}}|_{\sigma_{\mathfrak{a}}})^2 = (E^Y_{\mathfrak{a}}|_{\sigma_{\mathfrak{a}^*}})^2 =(\overline{E^Y_{\mathfrak{a}}}|_{\sigma_{\mathfrak{a}}})^2 = (\overline{E^Y_{\mathfrak{a}}}|_{\sigma_{\mathfrak{a}^*}})^2 = |E^Y_{\mathfrak{a}}|_{\sigma_{\mathfrak{a}}}|^2 = |E^Y_{\mathfrak{a}}|_{\sigma_{\mathfrak{a}^*}}|^2 &= (|E^*_{1,\psi}|^2)^Y.
\end{align*}

 Recall the explicit formula of $\varphi_{\mathfrak{aa}^*}$ and that $\psi$ is even, we can rewrite \eqref{DK5.1} as
\begin{align}\label{cleaned}
    \frac{12}{\Lambda^2(1,\psi)}\int_Y^{\infty}\int_0^1 y |E^*_{1,\psi}(z,\frac{1}{2})|^2 d\mu.
\end{align}

\subsection{Integral shift and poles of L-functions}
To calculate \eqref{cleaned}, recall
\begin{align*}
    E^*_{1,\psi}(z,\frac{1}{2}+iT) = 2\sqrt{y}\sum_{n\neq 0} \lambda_{1,\psi}(n,2iT) e(nx) K_{iT}(2\pi |n|y),
\end{align*} 
since $\psi(-1)=1$, the double integral in \eqref{cleaned} integral equals
\begin{align*}
    \int_Y^{\infty} 8y \sum_{n\geq 1} \lambda^2_{1,\psi}(n,0) K^2_{0}(2\pi ny) \frac{dy}{y}.
\end{align*}
With the integral and the sum interchanged, it equals
\begin{align}\label{periodint}
    8 (2\pi)^{-1} \sum_{n\geq 1} \frac{\lambda^2_{1,\psi}(n,0)}{n} g(2\pi nY),
\end{align}
where
\begin{align*}
    g(x):= \int_{x}^{\infty} y K_{0}^2(y) \frac{dy}{y}.
\end{align*}
The Mellin transform of $g$ is
\begin{align*}
    G(s) &= \int_0^{\infty} g(x)x^s \frac{dx}{x} = \int_0^{\infty} x^s \frac{dx}{x} \int_{x}^{\infty} y K_{0}^2(y) \frac{dy}{y} \\
    &= \int_{0}^{\infty} y K_{0}^2(y) \frac{dy}{y} \int_0^y x^s \frac{dx}{x} \\
    &= \frac{1}{s} \int_{0}^{\infty} y^{1+s} K_{0}^2(y) \frac{dy}{y},
\end{align*}
which according to the Mellin-Barnes formula [GR], 6.576.4, equals
\begin{align*}
    \frac{2^{-2+s}}{s \Gamma(1+s)} \Gamma
    ^4(\frac{1+s}{2}) .
\end{align*}
Applying the Mellin inverse of $G$, we can rewrite \eqref{periodint} by
\begin{align*}
    8 (2\pi)^{-1} \sum_{n\geq 1} \frac{\lambda^2_{1,\psi}(n,0)}{n}
    \frac{1}{2\pi i} \underset{(3)}{\int} G(s) (2\pi n Y)^{-s} ds
    = \pi^{-1}  \frac{1}{2\pi i} \underset{(3)}{\int}  \frac{\Gamma
    ^4(\frac{1+s}{2})}{s \Gamma(1+s)}  (\pi Y)^{-s} \sum_{n\geq 1} \frac{\lambda^2_{1,\psi}(n,0)}{n^{1+s}} ds.
\end{align*}
By Rankin-Selberg method of period integrals, the Dirichlet series equals
\begin{align*}
    \frac{\zeta(1+s) L^2(1+s, \psi) L(1+s, \chi_{0,N})}{L(2+2s,\chi_{0,N})} = \frac{\zeta^2(1+s) L^2(1+s,\psi)}{\zeta(2+2s)} \prod_{p\mid N}(1+\frac{1}{p^{1+s}})^{-1}.
\end{align*}

\noindent So, \eqref{cleaned} is equal to:
\begin{align}\label{3pole}
    \frac{12}{N \pi L^2(1,\psi)} \frac{1}{2\pi i} \underset{(3)}{\int}  \frac{\Gamma
    ^4(\frac{1+s}{2})}{s \Gamma(1+s)}  (\pi Y)^{-s} \frac{\zeta^2(1+s) L^2(1+s,\psi)}{\zeta(2+2s)} \prod_{p\mid N}(1+\frac{1}{p^{1+s}})^{-1} ds.
\end{align} 

\begin{remark}\label{explanation}
Now we can see why our case of $T=0$ and $\chi$ is quadratic is most complicated: the integrand has a triple pole at $s=0$, whereas in other cases it can have at most a double pole.  
\end{remark}

\noindent Since GRH implies $L(1+s,\psi) \ll \log\log N$ for $\Re s=0$, for each $N$ we have $c \in (0,\tfrac{\log\log\log N}{\log N})$ such that $L(1+s,\psi) = O(\sqrt{\log N})$ for $\Re s=-c$. Since the integrand of \eqref{3pole} has no poles other than $s=0$ for $\Re s \in (-c,3)$, we have
\begin{align*}
    \eqref{3pole} = \frac{12}{N \pi L^2(1,\psi)}  \Big( \underset{s=0}{\text{Res }} H(s) + \frac{1}{2\pi i} \underset{(-c)}{\int} H(s) ds \Big),
\end{align*}
with
\begin{align*}
    H(s)=\frac{\Gamma
    ^4(\frac{1+s}{2})}{s \Gamma(1+s)}  (\pi Y)^{-s} \frac{\zeta^2(1+s) L^2(1+s,\psi)}{\zeta(2+2s)} \prod_{p\mid N}(1+\frac{1}{p^{1+s}})^{-1}.
\end{align*}

\noindent For the integral we have
\begin{align*}
    \frac{12}{N \pi L^2(1,\psi)} \frac{1}{2\pi i} \underset{(-c)}{\int} H(s) ds \ll \frac{(\log\log N)^2}{N} Y^c L^2(1-c,\psi) \sum_{p\mid N} p^{-1+c} \ll \frac{\log N (\log\log N)^3}{N} Y^c.
\end{align*}

On the other hand, since
\begin{align*}
    \zeta(1+s) = \frac{1}{s} + \gamma_0 + \gamma_1 s + O(s^2),
\end{align*}
for some constant $\gamma_0,\gamma_1$, we have
\begin{align}\label{Laurentanal}
    \underset{s=0}{\text{Res }} H(s) = (2\gamma_1 + \gamma_0^2) K(0) + 2\gamma_0 K'(0) + K''(0),
\end{align}
where
\begin{align*}
    K(s) = \frac{\Gamma
    ^4(\frac{1+s}{2})}{\Gamma(1+s)}  (\pi Y)^{-s} \frac{ L^2(1+s,\psi)}{\zeta(2+2s)} \prod_{p\mid N}(1+\frac{1}{p^{1+s}})^{-1}.
\end{align*}

\noindent The following computation is straightforward (recall \eqref{A1}):
\begin{align*}
    K(0) &= 6 L^2(1,\psi) \frac{N}{\nu(N)} \\
    K'(0) &= K(0) \Big( -\log Y + 2\frac{L'}{L}(1,\psi) + \sum_{p\mid N}\frac{p^{-1} \log p}{1+p^{-1}} + O(1) \Big) \\
    K''(0) &= K(0) \Big( (\frac{K'(s)}{K(s)})'|_{s=0} + (\frac{K'(0)}{K(0)})^2 \Big) \\
    &= K(0) \Big( \sum_{p\mid N} \frac{p^{-1} \log^2 p}{1+p^{-1}} + O((\log\log N)^2) \Big).
\end{align*}
\noindent Thus \eqref{Laurentanal} yields
\begin{align*}
    \frac{12}{N \pi L^2(1,\psi)} \underset{s=0}{\text{Res }} H(s) = \frac{72}{\nu(N)} \Big( \sum_{p\mid N} \frac{p^{-1} \log^2 p}{1+p^{-1}} + O((\log\log N)^2) \Big).
\end{align*}
Since
\begin{align*}
    \sum_{p\mid N} \frac{p^{-1} \log^2 p}{1+p^{-1}} \leq \log N \sum_{p\mid N} \frac{\log p}{1+p} \ll \log N \log\log N,
\end{align*}
so \eqref{Laurentanal} is bounded by $\nu^{-1}(N) \log N \log\log N$. Furthermore,
\begin{align*}
    \eqref{3pole} \ll \frac{\log N (\log\log N)^3}{N} Y^c + \frac{\log N \log\log N}{\nu(N)} = o(\frac{\log^2 N}{\nu(N)}).
\end{align*}
This solves the first part of \eqref{broad}.
\subsection{Cauchy's inequality}
Now we are left to estimate
\begin{align*}
    \int_{\mathcal{F}} \Big( |E^Y_{\mathfrak{a}}|^2 \overline{E^Y_{\mathfrak{a}}} e_{\mathfrak{a}} + |E^Y_{\mathfrak{a}}|^2 E^Y_{\mathfrak{a}} \overline{e_{\mathfrak{a}}} \Big) d\mu = 2\int_{\mathcal{F}} (E^Y_{\mathfrak{a}})^3 e_{\mathfrak{a}} d\mu.
\end{align*}
By Cauchy's inequality we have
\begin{align*}
    \int_{\mathcal{F}} (E^Y_{\mathfrak{a}})^3 e_{\mathfrak{a}} d\mu \leq \Big( \int_{\mathcal{F}} (E^Y_{\mathfrak{a}})^2 e^2_{\mathfrak{a}} d\mu \Big)^{\frac{1}{2}} \Big( \int_{\mathcal{F}} (E^Y_{\mathfrak{a}})^4 d\mu \Big)^{\frac{1}{2}}.
\end{align*}
Our work on \eqref{cleaned} has shown the first integral is $o(\frac{\log^2 N}{\nu(N)})$, and the second integral has the same (big O) bound by assumption, so we have completed the proof.


\begin{thebibliography}{}
\bibitem{A} R. Askey, Beta integrals and the associated orthogonal polynomials, Lecture
Notes in Math., 1395, Springer-Verlag, Berlin, 1989, 84--121.

\bibitem{BKY} V. Blomer, R. Khan and M. Young, \emph{Distribution of mass of holomorphic cusp forms.}, Duke Math. J. 162(14): 2609--2644 (2013).
\bibitem{BM} V. Blomer and D. Milicevi\'c, \textit{The second moment of twisted modular $L$-functions.}, Geom. Funct. Anal. 25, 453--516 (2015).
\bibitem{BFKMM} V. Blomer, \'E. Fouvry, E. Kowalski, P. Michel, and D. Milićević. \emph{On moments of twisted L-functions.} Amer. J. Math., 139, no. 3 (2017), 707--768.
\bibitem{BLS}
A. Booker, M. Lee and A. Str\"ombergsson,
\emph{Twist-minimal trace formulas and the Selberg eigenvalue conjecture,}
 J. London Math. Soc., Volume 102, Issue 3, (2020) 1067--1134.
\bibitem{CFKRS} B. Conrey, D. Farmer, J. Keating, M. Rubinstein and N. Snaith, \emph{Integral moments of L-functions.}, Proc. London Math. Soc. (3) 91 (2005) 33--104.
\bibitem{DK1} G. Djankovi\'c and R. Khan,
\emph{A conjecture for the regularized fourth moment of Eisenstein series, J. Number Theory,} 182 (2018), 236--257.

\bibitem{DK2} G. Djankovi\'c and R. Khan,
\emph{On the Random Wave Conjecture for Eisenstein series,} Int. Math. Res. Not., 23 (2020), 9694--9716.
\bibitem{H} P. Humphries, \emph{Equidistribution in shrinking sets and $L^4$-norm bounds for automorphic forms}, Math. Ann. 371, 1497--1543 (2018).
\bibitem{HK} P. Humphries and R. Khan, \emph{On the Random Wave Conjecture for Dihedral Maaß Forms}, Geom. Funct. Anal., 30(1), 34--125.
\bibitem{I1}
H. Iwaniec,
\textit{Topics in Classical Automorphic Forms}, volume~17 of {\em Graduate Studies in Mathematics}. American Mathematical Society, Providence, RI, 1997.
\bibitem{I2}
H. Iwaniec,
\emph{Spectral Methods of Automorphic Forms (Second Edition)}, volume~53 of {\em Graduate Studies in Mathematics}. American Mathematical Society, Providence, RI, 2002.
\bibitem{PY} J. Pan and M. Young, \emph{Quantum Unique Ergodicity for Eisenstein Series in the Level Aspect.}, Commun. Math. Phys. 385, 227--266 (2021).
\bibitem{Sp} F. Spinu,
\emph{The $L^4$-norm of Eisenstein series,} PhD thesis, Princeton University.
\bibitem{Young} M. Young, \emph{Explicit calculations with Eisenstein series},
J. Number Theory,
199 (2019), 1--48.
\end{thebibliography}
\end{document}